\documentclass[10pt]{amsart}

\usepackage{amsmath,amsfonts,amsthm,amssymb}
\usepackage{color}
\usepackage[]{graphics}
\usepackage{verbatim}
\usepackage{eepic,epic}
\usepackage{epsfig,subfigure,epstopdf}
\usepackage{graphicx}

\numberwithin{equation}{section}

\newcommand{\al}{\alpha}

\def\Om{\Omega}

\def\Gd{\mathcal{G}_\delta}

\def\Ginf{\mathcal{G}_\infty}

\begin{document}

\title[Nonlocal-in-time model]
{Nonlocal-in-time dynamics and crossover of diffusive regimes}
\author[Q. Du  and Z. Zhou]{Qiang Du \and Zhi Zhou}
\begin{abstract}
We study a simple nonlocal-in-time dynamic system
 proposed for the effective modeling of complex diffusive
 regimes in  heterogeneous media.
We present its solutions and their
commonly studied statistics such as the mean square distance.
This interesting model employs a nonlocal operator to replace the
conventional first-order time-derivative. It introduces
a finite memory effect of a constant length encoded through a kernel function.
The nonlocal-in-time operator is related to fractional time derivatives that rely on the
entire time-history on one hand, while reduces to, on the other hand,
the classical time derivative if the length of the memory window diminishes.
This allows us to demonstrate the effectiveness of the nonlocal-in-time model
in capturing the crossover widely observed in nature
between the initial sub-diffusion and the long time normal diffusion.
\end{abstract}

\thanks{The research of the first author is supported in part by NSF DMS-1719699, NSF CCF-1704833,
  the ARO MURI Grant W911NF-15-1-0562 and the AFOSR  MURI center for Material Failure Prediction through peridynamics.
The research of the second author is also partially supported by the start-up fund provided by The Hong Kong Polytechnic University
and Hong Kong RGC grant No. 25300818.}

\maketitle

\section{Introduction}\label{sec:intro}
Diffusion in heterogeneous media
bears important implications in many  applications.
With the aid of single particle tracking,
recent studies have provided many examples of
anomalous diffusion, for example, sub-diffusion where the
spreading process happens
much more constricted and slower than the normal  diffusion
\cite{Berkowitz:2002,metzler2019brownian}.
Meanwhile,  the origins and mathematical
models of anomalous diffusion differ significantly
\cite{korabel2010paradoxes,sokolov2012models,berezhkovskii2014discriminating,mckinley2009transient}.
On one hand, new experimental standards
have been called for \cite{saxton2012wanted}. On the other hand, there are
needs for in-depth studies of mathematical models, many of
which are non-conventional and non-local \cite{sokolov2012models,Du:rev}.

Motivated by recent experimental reports on the crossover between initial transient
sub-diffusion and long time normal diffusion
in various settings \cite{he2016dynamic}, we present a simple dynamic equation that
provides an effective description of the diffusion
process encompassing these regimes. The main feature is to incorporate
memory effect or time correlations with a finite and fixed horizon length, denoted by $\delta>0$,  across the
dynamic process. The memory kernel is constant in space and  time so that  neither spatial inhomogeneities
nor time variations get introduced in the diffusivity coefficients.  This is  different from the approaches taken
in other models of anomalous diffusion like the variable-order fractional differential equation and
diffusing diffusivity \cite{chubynsky2014diffusing,grebenkov2019unifying,jain2017diffusing,Sun:2017}. In the model
under consideration here,
the memory effect dominates during the early time period, but as time goes on, the
fixed memory span becomes less significant over the long life history. As a result,
the transition from sub-diffusion to normal diffusion occurs naturally. A rigorous demonstration of this intuitive picture 
will be given later in more details.
We note that while the nonlocal-in-time diffusion equation may be related to
fractional diffusion equations \cite{Metzler:2000,metzler2004restaurant,sokolov2012models}
by taking special memory kernels \cite{Allen:2015},  they in general provide a new class of
models that effectively serve as a bridge between  anomalous diffusion and
normal diffusion, with the latter being a limiting case as the horizon length $\delta\to 0$.

Specifically,  let $\Delta$ be the
Laplacian (diffusion) operator in the spatial variable $x$ and
 $g=g(x,t)$ represent the initial (historical) data,
 we consider the following nonlocal-in-time diffusion equation
 for $u=u(x, t)$:
\begin{eqnarray}\label{eqn:fde}
   \Gd u  -\Delta u = 0,&  \qquad t>0, \nonumber\\
    u=g,& \qquad  t\in (-\delta,0).
\end{eqnarray}
The nonlocal
 operator $\Gd$ in (\ref{eqn:fde}) is defined by
\begin{equation}\label{eqn:op}
    \Gd v (t)=
     \int_0^\delta  \frac{v(t)-v(t-s)}{s}\rho_\delta(s)\,ds,
 \end{equation}
where the memory kernel function
$\rho_\delta=\rho_\delta(s)$  is assumed to be nonnegative with a compact support in $(0, \delta)$
and is integrable in $(0,\delta)$. 
In case that $s^{-1}\rho_\delta(s)$ is unbounded only at the origin,
the integral in (\ref{eqn:op}) should be interpreted as the limit of the integral of the same
 integrand over $(\epsilon, \delta)$ for $\epsilon>0$
 as $\epsilon \to 0$, where such a limit exists in an appropriate mathematical sense.

The nonlocal operator $\Gd$ \cite{DuYangZhou:2016,DuTaoTianYang:2015a, Du:stochastic}
forms part of the nonlocal vector calculus \cite{DuGLZ:2013,DuMengesha:2015}.
The positive nonlocal horizon parameter $\delta$ appearing in (\ref{eqn:op})
represents the range of nonlocal interactions or memory span.
For suitably chosen kernels, as $\delta\rightarrow0$,
nonlocal and memory effects diminish, so that the zero-horizon limit of the
nonlocal operator $\Gd u$ corresponds to the standard first order derivative $\frac{d}{dt}u$.
In particular, under the normalization condition
\begin{equation}\label{eqn:2moment}
    \int_0^\delta  \rho_\delta (s)\,ds  = 1,
\end{equation}
the kernel function $\rho_\delta (s)$ can be seen as a probabilistic density function(PDF)  defined for $s$ in $(0,\delta)$,
so that the nonlocal operator (\ref{eqn:op}) can be viewed as a "continuum" average of the backward difference
operators over the memory span measured by the horizon $\delta>0$.
The local limit, i.e., the standard derivative, is simply the extreme case
where $\rho_\delta(s)$ degenerates into a singular point measure at $s=0$.
In fact, one can see from a formal Taylor expansion that
\begin{eqnarray*}
\Gd  v(t) &= &  \frac{dv}{dt}(t)
+ \sum_{k=2}^\infty \frac{v^{(k)} (t)}{k!}
\int_0^\delta
{(-s)^{k-1}}
\rho_\delta(s)\,ds\\
&=& \frac{dv}{dt}(t) + O(\delta)\to \frac{dv}{dt}(t), \quad\mbox{as } \,\delta\to 0.
\end{eqnarray*}
Then the nonlocal-in-time model (\ref{eqn:fde}) recovers the classical (local) diffusion model.
\begin{equation*}
    \partial_t  u - \Delta u = 0.
\end{equation*}

Meanwhile, in another extreme case that $\delta\rightarrow\infty$, we may let the initial (historical) data $v(t)=v(0)$
for all $t\in(-\infty,0)$, and  take the kernel function to be of the fractional type,
i.e., 
\begin{equation}
\label{fractype}
\rho_\delta(s) = \frac{\al}{\Gamma(1-\al)} {s^{-\al}} \chi_{(0,\delta)}(s), \;\mbox{ for some }\; \al\in(0,1),
\end{equation}
where $\chi=\chi_{(0,\delta)}(s)$ denotes  the characteristic (indicator) function of $(0,\delta)$. Then,
the nonlocal operator $\Gd$ reproduces the Marchaud fractional derivative  of order $\alpha$,
which is equivalent to the Caputo fractional derivative for smooth functions at time $t>0$ (e.g. \cite[p.\,91]{KilbasSrivastavaTrujillo:2006}).
As a consequence, (\ref{eqn:fde}) recovers the fractional sub-diffusion model for $u=u(x,t)$:
\begin{eqnarray*}
   \partial_t^\alpha u  -\Delta u = 0,&  \qquad t>0,\\
    u(0)=g(x).
\end{eqnarray*}
The fractional subdiffusion model has often been used to describe the continuous time random
walk (CTRW) of particles in heterogeneous media, where trapping events occur.
In particular,  particles get repeatedly immobilized in the
environment for a trapping time drawn from the waiting time probability density function has a heavy tail,
i.e., $\omega(t)\propto t^{-\alpha-1}$ as $t\rightarrow\infty$.
\cite{Metzler:2000,metzler2004restaurant}.
As discussed later in Section \ref{sec:3},  
the nonlocal-in-time model (\ref{eqn:fde}) can be related to a trapping model
where the kernel function
describes the distribution of waiting time probability.



In comparison with the classical diffusion equation
(the $\delta\to 0$ local limit) and fractional diffusion (the $\delta\to \infty$ limit),
the nonlocal-in-time evolution equation 
provides a more general model and an interesting intermediate case
 to study the "finite history dependence" with a given $\delta\in(0,\infty)$.
The goal of this paper is to present  the modeling capability and the behavior
of the solutions of the nonlocal-in-time dynamics (\ref{eqn:fde}),
so as to demonstrate how the simple PDE model can serve as a bridge linking
the standard diffusion with the fractional sub-diffusion.
We refer to \cite{DuYangZhou:2016,DuTaoTianYang:2015a,DuMengesha:2015,DTZ2019} for
more development on the mathematical background and numerical analysis
of the nonlocal operators and the nonlocal-in-time dynamic systems.

\section{Solutions of the nonlocal-in-time model}\label{sec:2}
To begin with, let us study the solution behavior of the
nonlocal-in-time diffusion dynamics. To correlate with data observed from
physical and biological experiments, we focus on  the spatial solution
profiles at various stages and the time evolution of the mean square displacement. In order to
make comparisons with normal and fractional diffusion models, we choose
fractional type memory kernels of the type (\ref{fractype}) as illustrations.

\subsection{Fundamental solutions}
We start with the model (\ref{eqn:fde}) defined spatially on the
real line, with a time stationary Dirac-delta measure in $x$ as the initial distribution, i.e.,
$g(x,t)=\delta(x)$ for $x\in (-\infty,\infty)$ and $ t\in(-\delta,0)$.
The solution to (\ref{eqn:fde})  in this case may be called the fundamental solution
of the nonlocal model. Indeed, we can take the kernel function to have a unit integral,
as  assumed in Section I.
Then, the limiting local solution is the well-known fundamental solution of the heat equation, expressed by a Gaussian function
\begin{equation*}
    u_0(x,t) = \frac{1}{\sqrt{4\pi t} } e^{-x^2/4t}.
\end{equation*}

As a comparison, in another limit, as $\delta$ approaches $\infty$, with the kernel function
\begin{equation*}\rho_\delta(s) = \frac{\al}{\Gamma(1-\al)} {s^{-\al }}\chi_{(0,\delta)}(s),
\end{equation*}
the nonlocal derivative reproduces  a fractional derivative, i.e., $\Gd\rightarrow \Ginf = \partial_t^\alpha$, the fundamental solution is given
by the Fox H function
\begin{equation*}
    u_\infty(x,t) = \frac12 t^{-\alpha/2} P_{\alpha}(x/t^{\alpha/2}),
\end{equation*}
with $P_{\alpha}(y) \sim A y^a e^{-by^c}$ as $ y\rightarrow \infty$ for
    constants  $a=\frac{2\alpha-2}{2-\alpha}$, $b = (2 -\alpha) 2^{-\frac{2}{2-\al}}\al^{\frac{\al}{2-\alpha}}$,
    and $c=\frac2{2-\al}$ \cite{mainardi2001fundamental}.
  
%
For the nonlocal-in-time model (\ref{eqn:fde}), let $\widetilde u(\xi,t)$ denote the
Fourier transform of fundamental solution $u(x,t)$ with respect to $x$,  we get a scalar nonlocal-in-time initial value problem
$$ \Gd \widetilde u(\xi,t) + \xi^2 \widetilde u (\xi,t) = 0$$ for $t>0$ with an initial data
$\widetilde u (\xi,t) = 1$ for $t\in(-\delta, 0)$.
The Laplace transform of $\widetilde u(\xi,t)$ with respect to time $t$ is given by
$$   {\widehat {\widetilde u}} (\xi,z) = \frac{z^{-1}K(z)}{K(z)+\xi^2},\;\mbox{with}\;
 K(z) = \int_0^\delta (1-e^{-zs})s\rho_\delta(s)\,ds.
 $$
Applying the inverse Laplace and Fourier transforms,
we obtain a formal analytic representation of the solution
\begin{eqnarray}
    u(x,t) &=& \frac{1}{2\pi}  \int_{-\infty}^{\infty} e^{\rm{i} \xi x } \frac{1}{2\pi\rm{i}}
    \int_{\sigma-\rm{i}\infty}^{\sigma+\rm{i}\infty}
    \frac{z^{-1}K(z)}{K(z)+\xi^2} e^{zt}\,dz  \,d\xi \nonumber\\
   &=& \frac{1}{2\pi\rm{i}}\int_{\sigma-\rm{i}\infty}^{\sigma+\rm{i}\infty} e^{zt} \frac{1}{2\pi} \int_{-\infty}^{\infty} e^{\rm{i} \xi x }
   \frac{z^{-1}K(z)}{K(z)+\xi^2}\,d\xi \,dz\nonumber\\
    & =&
     \frac{1}{2\pi\rm{i}}\int_{\sigma-\rm{i}\infty}^{\sigma+\rm{i}\infty} e^{zt} \widehat u(x,z) \,dz.
     \label{eqn:fund}
\end{eqnarray}
The  above integrals  are well-defined for $(x,t)\in \mathbb{R} \times(0,\infty)$.
The fundamental solution $u=u(x,t)$ is continuous in $x$ and $t$ and
piecewise smooth in $x$ away from $x=0$ for $t>0$.
For $z>0$, we can get
$$ \widehat u(x,z) = \frac{ \sqrt{K(z)}}{2z} e^{-|x| \sqrt{K(z)}}.$$

\begin{figure}[hbt!]
\centering
\includegraphics[trim = .1cm .1cm .1cm .1cm, clip=true,width=0.55\textwidth]{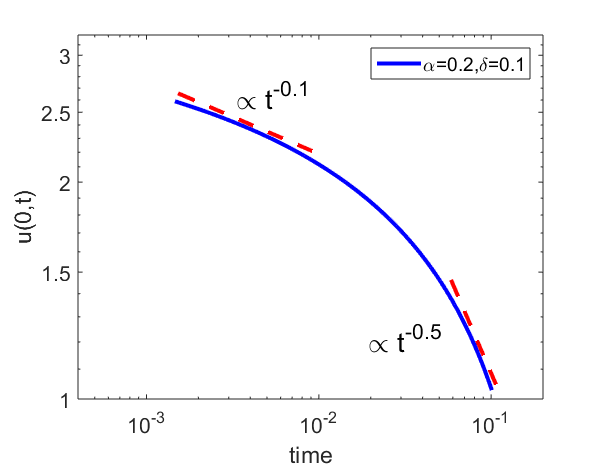}
\vspace{-0.1cm}
 \caption{$u(0,t)$ of nonlocal-in-time diffusion with $\delta=0.1$, $\alpha=0.2$.}\label{fig:u0t}
\end{figure}

Concerning the time variation of the fundamental solution of (\ref{eqn:fde}) with the kernel (\ref{eqn:rho}), for example at $x=0$, we note first
$ \widehat u(0,z)$ is monotone decreasing in $z$ and
 \begin{equation*}
      \widehat u(0,z)
  \approx	
  \left\{\begin{array}{lr}
	\sqrt\frac{\Gamma(2-\alpha)}{4\alpha\delta^{1-\alpha}} z^{\alpha/2-1}
	  &\quad \mbox{as}~~z\rightarrow \infty,  \\
     \frac12z^{-1/2}
  &\quad \mbox{as}~~z\rightarrow 0.
 \end{array}\right.
\end{equation*}
 Then the Karamata-Feller-Tauberian theorem \cite{Feller:2008}
gives
\begin{equation*}
   u(0,t) \approx	
  \left\{\begin{array}{lr}
\frac{\sqrt{\Gamma(2-\alpha)}}{\sqrt{4\alpha\delta^{1-\alpha}}\Gamma(1-\alpha/2)} t^{-\alpha/2}
	  &\quad \mbox{as}~~t\rightarrow 0,
  \\
   \quad  \frac1{\sqrt{4\pi t}}
  &\quad \mbox{as}~~t\rightarrow \infty.
 \end{array}\right.
\end{equation*}
Again, as $t\rightarrow \infty$, the decay property is the same as the standard heat kernel, while
 as $t\rightarrow 0$, the solution behaves like $O(t^{-\alpha/2})$,
which is the same as the fundamental solution of fractional sub-diffusion.
As an illustration, Fig~\ref{fig:u0t} shows a plot of numerical solution $u(0,t)$ for $\alpha=0.2$
and $\delta=0.1$.

\subsection{Mean square displacement (MSD)}
An interesting characteristic of the diffusion process is the mean square displacement  (MSD), denoted by $m(t)$ and
defined by
\begin{equation*}
  m(t) := \int_{\mathbb{R}} x^2 u(x,t) \,dx.
\end{equation*}
For the fundamental solution (\ref{eqn:fund}) of (\ref{eqn:fde}) with
the kernel $\rho_\delta$ being a normalized PDF, i.e., (\ref{eqn:2moment}) is
satisfied,
the corresponding  $m(t)$
satisfies the nonlocal initial value problem
$\Gd m(t) =2$ for $t>0$ with $m(t) = 0$ for $t\in (-\delta, 0)$. One can get that $m(t)\approx 2t$ as
$t\to \infty$. Meanwhile, we have
\begin{equation*}
 m(t)  \approx	  \frac{\sin(\alpha\pi)\delta^{1-\alpha}}{(1-\alpha)\pi}t^{\alpha}\, \quad \mbox{as}\;~ t\to 0
\end{equation*}
if we use special kernels of the fractional type:
\begin{equation}\label{eqn:rho}
  \rho_\delta(s)= (1-\al) \delta^{\al-1} {s^{-\al }}.
\end{equation}
Thus, we observe the transition from the sub-diffusion initially to normal diffusion at a later time.

\begin{figure}[hbt!]
\centering
\includegraphics[trim = .1cm .1cm .1cm .1cm, clip=true,width=0.50\textwidth]{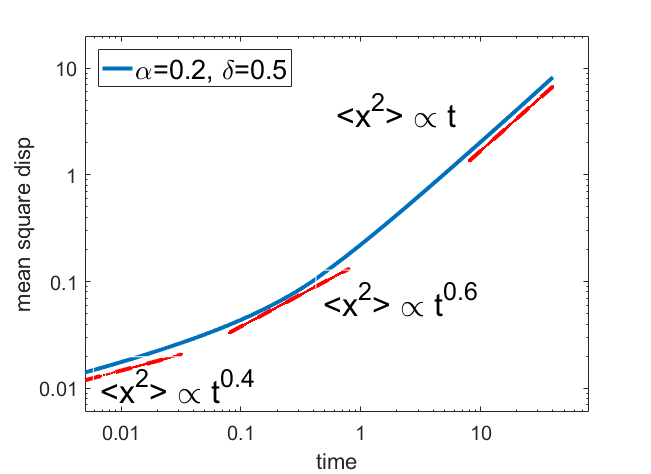}
\includegraphics[trim = .1cm .1cm .1cm .1cm, clip=true,width=0.45\textwidth]{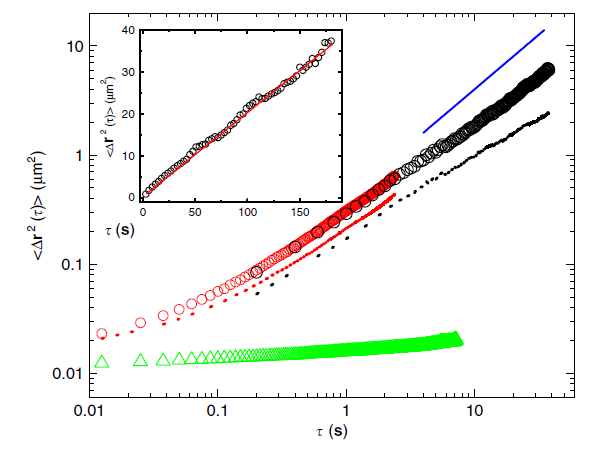}
\vspace{-0.1cm}
 \caption{(a) The plot of MSD of nonlocal-in-time diffusion model (\ref{eqn:fde}) with $\alpha=0.2$ and $\delta=0.5$;
 (b) The immobile AChRs: Crossover from sub-diffusion to normal diffusion observed from the MSD curve \cite[Fig. 3]{he2016dynamic}.}\label{fig:msd}
\end{figure}

In figure \ref{fig:msd}, we plot the numerical solution of $\Gd m(t) = 2$, i.e., the mean square displacement of the nonlocal model corresponding
to $\rho_\delta(s)$ given by (\ref{eqn:rho})
with $\alpha=0.2$ and $\delta=0.5$.
The experiment again illustrates the analytically suggested change
from the early fractional anomalous diffusion regime to
the later standard diffusion regime.
This  "transition" or "crossover" behavior appears in many practical applications, e.g.
diffusion in
lipid bilayer systems of varying chemical compositions \cite[Fig.~2]{jeon2012anomalous},
and lateral motion of the acetylcholine receptors on live muscle cell membranes\cite[Figs.~3,~4]{he2016dynamic}.

\begin{figure}[hbt!]
\centering
{
\includegraphics[trim = .1cm .1cm .1cm .1cm, clip=true,width=0.31\textwidth]{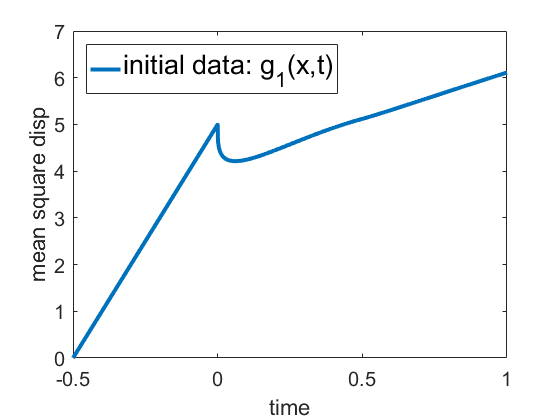}}
{
\includegraphics[trim = .1cm .1cm .1cm .1cm, clip=true,width=0.31\textwidth]{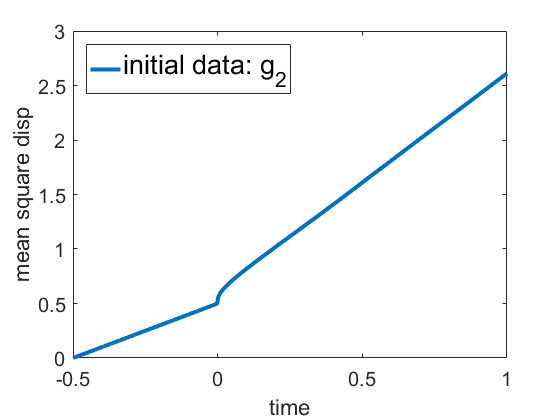}}
{
\includegraphics[trim = .1cm .1cm .1cm .1cm, clip=true,width=0.31\textwidth]{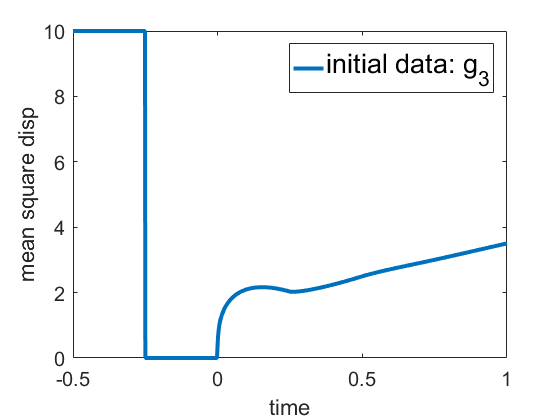}}
\vspace{-0.1cm}
 \caption{Mean square displacement of nonlocal-in-time model with different initial distribution, $\alpha=0.2$, $\delta=0.5$}\label{fig:ini-msd}
\end{figure}

We can further study the dependence of MSD on the initial  historical data.
Based on the earlier discussion,
we know that if the initial condition $g$ is history-independent, i.e., $g(x,t):=u_0(x)$,
then $m(t)$, the mean square displacement of  the nonlocal-in-time diffusion with a kernel function (\ref{eqn:rho}),
is  increasing and exhibits a weak singularity at the start, as its fractional counterpart.
However, the nonlocal-in-time model (\ref{eqn:fde})
allows a history-dependent initial data,
which will affect the growth behavior of the mean square displacement.
In Fig. \ref{fig:ini-msd}, we plot the solution of initial value problem $\Gd m(t)=2$, with different initial (historical) data for $t\in(-\delta,0)$:
\begin{equation*} 
g_1(t) = 5(1+2t),\quad g_2(t) = \frac{1+2t}{2}  \quad \mbox{and} \quad g_2(t)=10*\chi_{[-\delta,\delta/2]}.
\end{equation*}

Numerical results show that the solution $m(t)$ of initial data $g_1$
decays in the beginning, due to the historical-dependence of the nonlocal-in-time dynamics (where the initial data grows rapidly),
and then increases at a later time,
while for the initial data $g_2$, whose slope is smaller, we observe a strictly increasing mean square displacement function.
In fact, one can prove that for initial data $k(1+2t)$, the solution $m(t)$ will keep strictly increasing when $k\le 1$. 
Besides, for the step initial data $g_3$, the solution increases dramatically near $t=0$,
and then decreases a little before its constantly linear growth. The growth behavior of $m(t)$ for
different initial data is interesting and awaits further theoretical study.

\subsection{More on the long time normal diffusion limit}
The long-time normal diffusion behavior can also be observed directly from the model (\ref{eqn:fde}) by a simple scaling.
In particular, for a typical rescaled density $\rho_\delta(s)=\delta^{-1} \rho_1(\delta^{-1}s)$ with
$\rho_1$ being a density over the unit interval, then for $\delta_0$ and
$v(x,t)=u(x'.t')$ with $x'=x\delta_0^{-1/2}$ and
$t'=t \delta_0^{-1}$, it holds that
\begin{eqnarray*}
    \Gd  v (x, t)&=&
     \int_0^\delta  \frac{v(x,t)-v(x,t-s)}{s}\rho_\delta(s)\,ds\\
     &=& \frac{1}{\delta}
     \int_0^\delta  \frac{v(x,t)-v(x,t-s)}{s}\rho_1(\frac{s}{\delta})\,ds\\
           &=&  \frac{1}{\delta}   \int_0^\delta  \frac{u(\frac{x}{\sqrt{\delta_0}},\frac{t}{\delta_0})-u(\frac{x}{\sqrt{\delta_0}},\frac{t-s}{\delta_0})}{s}\rho_1(\frac{s}{\delta})\,ds\\
      &=&   \frac{1}{\delta_0}   \int_0^{\delta'}  \frac{u(x',t')-u(x',t'-s')}{s'}\rho_{\delta'} (s')\,ds'\\
      &=&   \frac{1}{\delta_0}  \mathcal{G}_{\delta'} u(x', t'),
\end{eqnarray*}
 by a change of variables $\delta=\delta' \delta_0$ and $s=s'\delta_0$.
Then, using  $\Delta_{x} v(x,t) = \delta_0^{-1} \Delta_{x'}v(x',t')$,
   we see that for any solution $u=u(x,t)$ of  (\ref{eqn:fde}),
   the rescaled solution $v$ satisfies the same equation in the rescaled variables corresponding to a kernel with a rescaled horizon
   $\delta'=\delta/\delta_0$.
   If we let $\delta_0\to \infty$, then
   the rescaled horizon $\delta'$ goes to zero, so that $\mathcal{G}_{\delta'}$ is effectively the conventional local time derivative.
  We thus  see that  with a diminishing memory effect, $v$ approximately satisfies the classical normal diffusion on the $O(t)$ scale, so
   does $u$ on the long time scale $O(t')>>O(t)$.

For a fixed nonlocal horizon $\delta$ and at a given time $t$, one may derive the asymptotic behavior of the
solution to the nonlocal-in-time model (\ref{eqn:fde}), as $x\rightarrow\infty$, from its Laplace transform.
In particular, we have
\begin{equation*}  \widehat u(x,z) \approx \frac12 c_\delta z^{\alpha/2-1} e^{-c_\delta|x|z^{\alpha/2}}, \;\mbox{ as }\;
z\rightarrow\infty,\end{equation*}
for $c_\delta=\sqrt{\Gamma(2-\alpha)\delta^{\alpha-1} \alpha^{-1}}$. Hence by inverse Laplace transform, we have
\begin{equation*} u(x,t) \approx \frac12 c_\delta t^{-\alpha/2} P_{\alpha/2}(c_\delta x/t^{\alpha/2}),\;\mbox{ as }\; t\rightarrow 0,\end{equation*}
with the Fox H function $P_\alpha$.
This implies that  $u(x,t)$ has a fractional exponential tail, i.e., for a fixed $t>0$,
\begin{equation*} u(x,t)\approx A_{\alpha,\delta,t} e^{-b_{\alpha,\delta,t}x^{2/(2-\alpha)}},\; \mbox{ as  }\; x\rightarrow\infty.\end{equation*}

\subsection{Smoothing properties}
Besides the patterns on statistics like MSD discussed above,
another interesting feature of the nonlocal-in-time dynamics is its gradual smoothing property,
namely, the solutions can become more and more smooth (as functions of the spatial variables) 
as time goes on. This can be derived similarly as 
 in \cite{DuYangZhou:2016}.
In fact, with the fractional kernel function in (\ref{eqn:rho}), we may see  that
\begin{equation*} | \widetilde u(\xi,t) | \le {c}/(1+b_\delta |\xi|^2 t^\alpha)\end{equation*}
for any $t>0$ and some constants $c$ and $b_\delta$
\cite[Theorem 3.2]{DuYangZhou:2016}.
This implies  that $|\xi|^s | \widetilde u(\xi,t) |$ is square-integrable in
$(-\infty,\infty)$ only if $s<\frac32$.
This restriction reflect the limiting smoothness of the fundamental solution,
i.e., the solution $u(x, t)$, for any $t>0$, has nearly $3/2$-order square integrable fractional
spatial derivative  and hence is only piecewise smooth.
In fact,  it roughly gains two more orders of differentiability
with each additional $\delta$ increment in time.
This observation indicates that  the smoothing of the solution $u(x,t)$  takes place incrementally over time.

\subsection{Comparison with fractional and normal diffusion via numerical illustrations}

\begin{figure}[hbt!]
\centering
\subfigure
{\includegraphics[trim = .1cm .1cm .1cm .1cm, clip=true,width=0.32\textwidth]{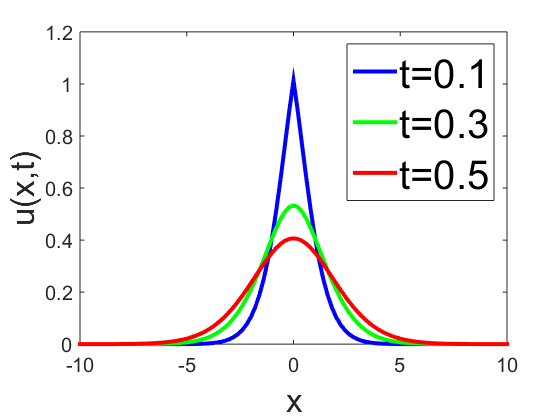}}
\subfigure
{\includegraphics[trim = .1cm .1cm .1cm .1cm, clip=true,width=0.32\textwidth]{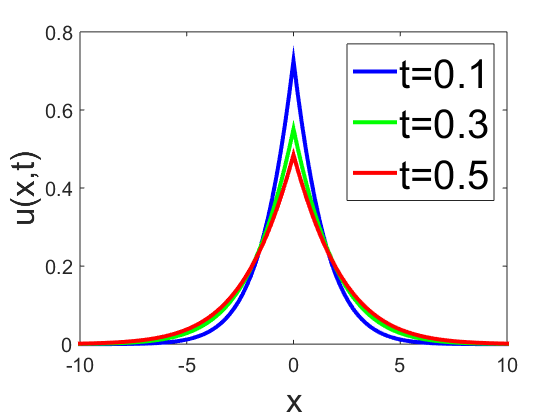}}
\subfigure
{\includegraphics[trim = .1cm .1cm .1cm .1cm, clip=true,width=0.32\textwidth]{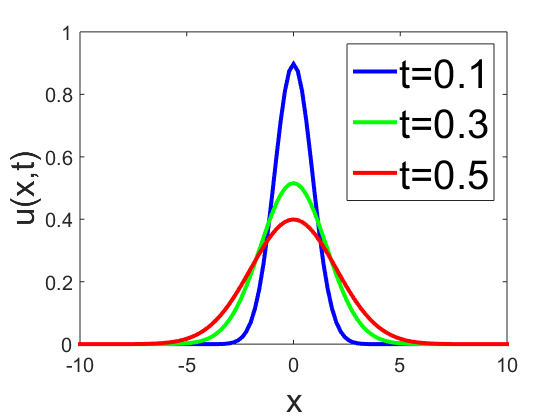}}
\vspace{-0.1cm}
 \caption{Left: Numerical solutions of the nonlocal-in-time  diffusion equation
  with $\delta=0.2$, $\alpha=0.5$.
 Middle: Numerical solutions of the fractional diffusion with $\alpha=0.5$.
 Right: Numerical solutions of the local diffusion.}\label{fig:inf_time}
\end{figure}

In Fig.~\ref{fig:inf_time}, we plot the numerical solution of the nonlocal-in-time model
with $\delta=0.2$, $\alpha=0.5$ at different time, and compare with those
of the fractional diffusion and normal diffusion.

From numerical experiments, we observe that the solution is piecewise smooth spatially for $t\in(0,\delta)$ (the blue curve in Fig.~\ref{fig:inf_time}(a)),
and is  getting more regular as $t$ keeps increasing (the green and red curves in Fig.~\ref{fig:inf_time}(a)).
This marks another distinct feature of the nonlocal-in-time
 model (\ref{eqn:fde})  with a finite memory that it has the same spatial smoothing property
as the corresponding anomalous sub-diffusion (Fig.~\ref{fig:inf_time}(b)), initially,
but gets improved smoothing incrementally in time and
exhibits the smoothing behavior of standard diffusion (Fig.~\ref{fig:inf_time}(c)) as $t$ approaches infinity.

In making the above comparisons, we note that the kernel $\rho_\delta$
in (\ref{eqn:rho}) differs from the $\rho_\infty$ used for the fractional diffusion
by a constant factor. Similar comparisons can be made for a rescaled kernel,
and it is easy to see that, if a new kernel is taken as $\sigma \rho_\delta$
with a constant factor $\sigma>0$,
then the fundamental solution of the nonlocal model with the new kernel is simply
given by $u(x/\sqrt{\sigma}, t)$
with $u(x,t)$ being  the fundamental solution of (\ref{eqn:fde}) with the original kernel.

\begin{figure}[hbt!]
\centering
\subfigure[t=0.1]
{\includegraphics[trim = .1cm .1cm .1cm .1cm, clip=true,width=0.45\textwidth]{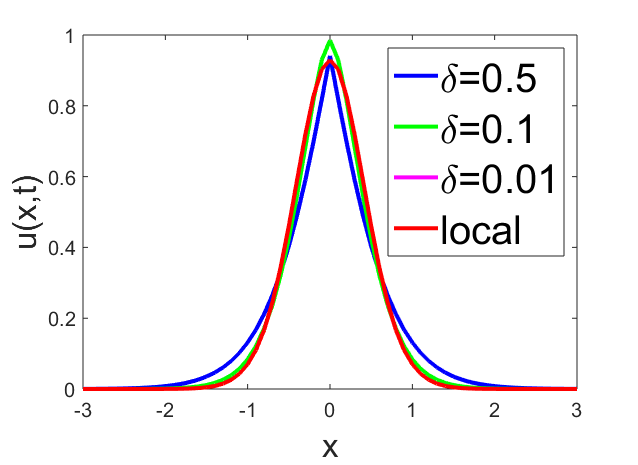}}
\subfigure[t=0.5]{
\includegraphics[trim = .1cm .1cm .1cm .1cm, clip=true,width=0.45\textwidth]{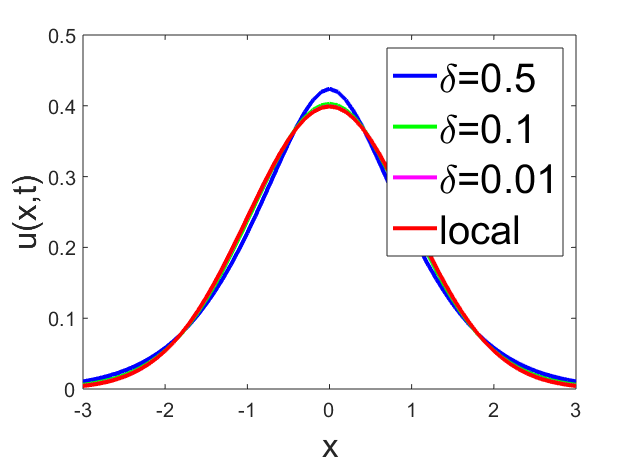}}
\subfigure[t=0.1]
{\includegraphics[trim = .1cm .1cm .1cm .1cm, clip=true,width=0.45\textwidth]{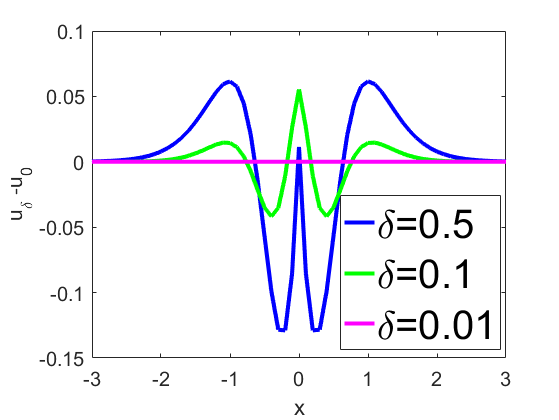}}
\subfigure[t=0.5]{
\includegraphics[trim = .1cm .1cm .1cm .1cm, clip=true,width=0.45\textwidth]{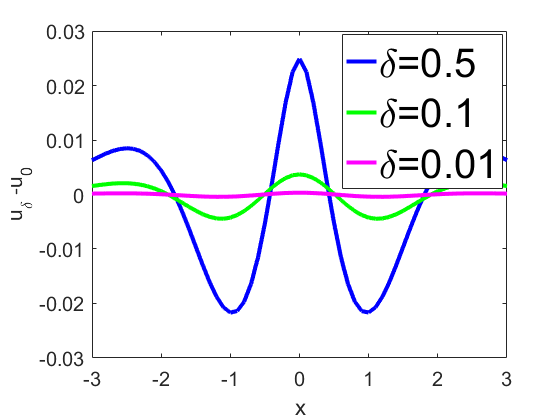}}
\vspace{-0.1cm}
 \caption{Figures (a) and (b): The plot of numerical solutions of nonlocal-in-time diffusion at $t=0.1$, $0.5$ with
 $\alpha=0.5$ and different nonlocal horizons.
 Figures (c) and (d): The plot of difference between solutions of nonlocal diffusion with different nonlocal horizons
 and the solution of local diffusion.}\label{fig:local_lim_inf}
\end{figure}

It is interesting to numerically study the local limit of the nonlocal-in-time model that has been
rigorously established in \cite{DuYangZhou:2016}.
In Fig. \ref{fig:local_lim_inf}, we plot the numerical solutions of nonlocal-in-time model
with different nonlocal horizons at $t=0.1$
and $0.5$, as well as the difference between the nonlocal solutions and the local one.  
It can be observed that as $\delta$ goes to zero, the solution of the nonlocal diffusion model
converges to the solution of the local one.

\subsection{Nonlocal-in-time dynamics on a finite spatial domain}
To complement the study on the infinite one dimensional spatial domain,
 we now briefly turn to the nonlocal-in-time parabolic equation (\ref{eqn:fde})
in a bounded domain with a homogeneous Dirichlet boundary condition.
Specifically, we consider the following initial-boundary value problem for $u=u(x, t)$:
\begin{eqnarray}\label{eqn:fde2}
   \Gd u(x,t)- \Delta u (x,t)= 0,&\quad \mbox{in} ~~\Om\times[0,T],\nonumber\\
    u(x,t)=0, &\quad \mbox{in} ~~\partial\Om\times[0,T], \\
    u(x,t)=g(x,t), &\quad \mbox{in}  ~~\Om\times[-\delta,0),\nonumber
\end{eqnarray}
where $\Om$ is a bounded convex polyhedral domain in $\mathbb{R}^d$ $(d \geq 1)$ with a boundary 
$\partial\Om$. Mathematical studies of (\ref{eqn:fde2}) such as well-posedness and numerical analysis can be found in \cite{DuYangZhou:2016}.

\begin{figure}[hbt!]
\centering
{
\includegraphics[trim = .1cm .1cm .1cm .1cm, clip=true,width=0.32\textwidth]{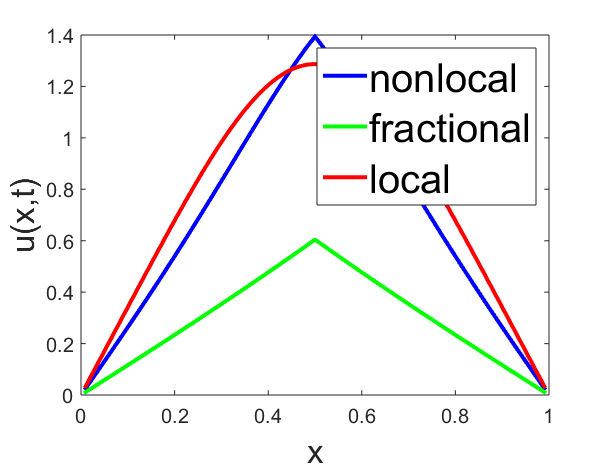}}
{
\includegraphics[trim = .1cm .1cm .1cm .1cm, clip=true,width=0.32\textwidth]{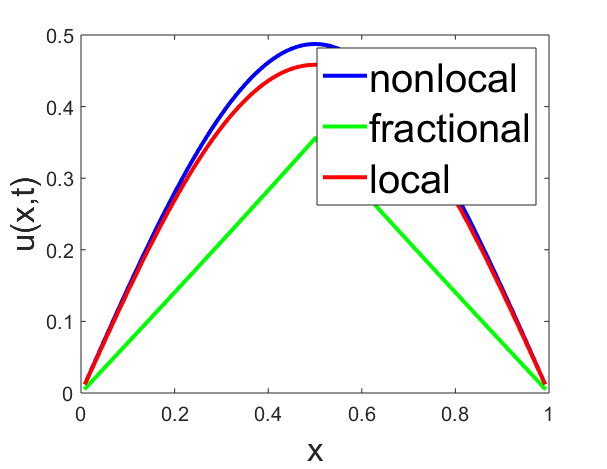}}
{
\includegraphics[trim = .1cm .1cm .1cm .1cm, clip=true,width=0.32\textwidth]{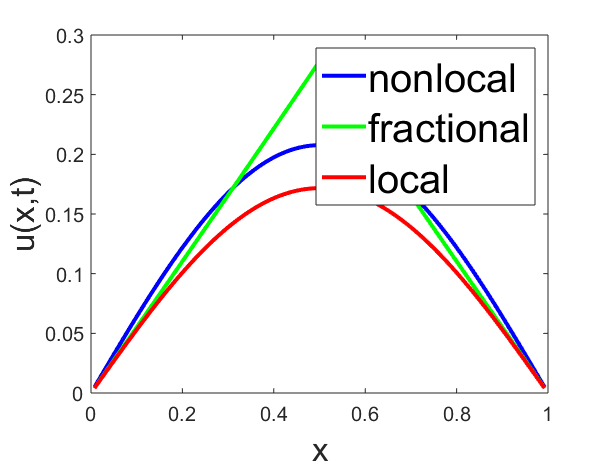}}
\vspace{-0.1cm}
 \caption{Numerical solutions of local, fractional and nonlocal diffusion models,  $\delta=0.1$ and $\alpha=0.5$.}\label{fig:case2.1}
\end{figure}

To illustrate our findings, we use the fractional kernel given by (\ref{eqn:rho})
and compare its solution behavior with those of fractional diffusion and
local diffusion.
In Fig.~\ref{fig:case2.1}, we present  a numerical
solution of an 1-D nonlocal model (\ref{eqn:fde2})  on   $\Omega=(0,1)$
at different time  $t=0.05$, $0.15$ and $0.25$  with  $\delta=0.1$ and $\alpha=0.5$, The initial (historical) data
  is taken as a Dirac-delta measure at $x=1/2$.
We can observe that the nonlocal diffusion gradually regularize the solution as $t\rightarrow\infty$.
Moreover, the fractional diffusion decays fastest for small $t$, and the classical
local diffusion decays fastest for large $t$, while nonlocal diffusion exhibits an intermediate behavior.

\begin{figure}[h!]
\subfigure[ {\scriptsize nonlocal, $t=0.01$}]{
  \includegraphics[trim = .1cm .1cm .1cm .1cm, clip=true,width=0.3\textwidth]{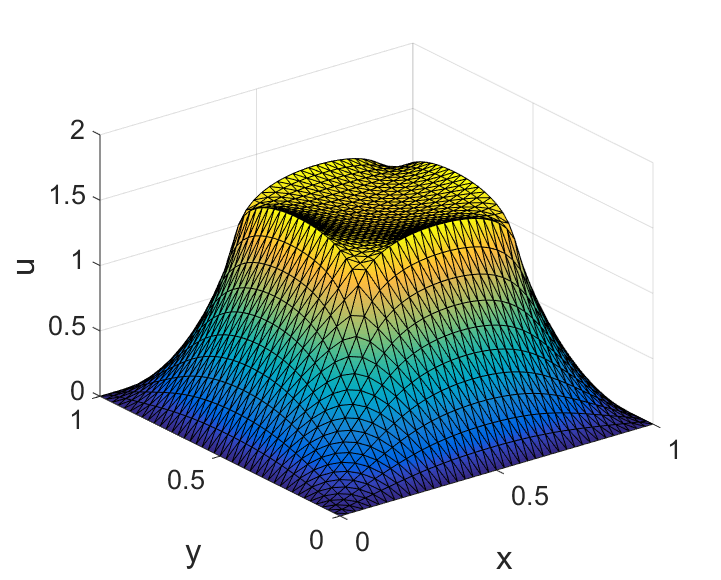}}
\subfigure[{\scriptsize fractional, $t=0.01$}]{
  \includegraphics[trim = .1cm .1cm .1cm .1cm, clip=true,width=0.3\textwidth]{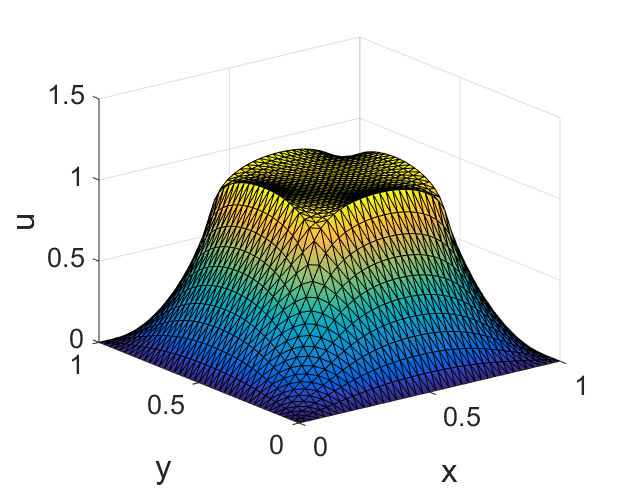}}
\subfigure[{\scriptsize standard, $t=0.01$}]{
  \includegraphics[trim = .1cm .1cm .1cm .1cm, clip=true,width=0.3\textwidth]{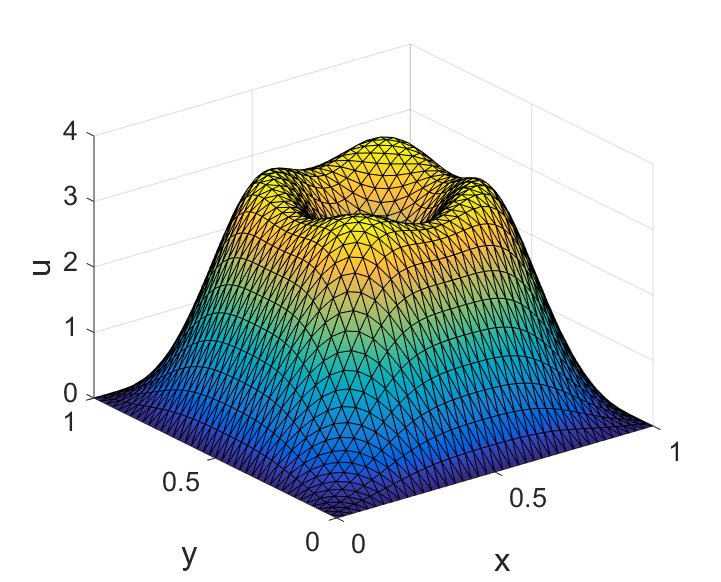}}
\subfigure[{\scriptsize nonlocal, $t=1.1$}]{
  \includegraphics[trim = .1cm .1cm .1cm .1cm, clip=true,width=0.3\textwidth]{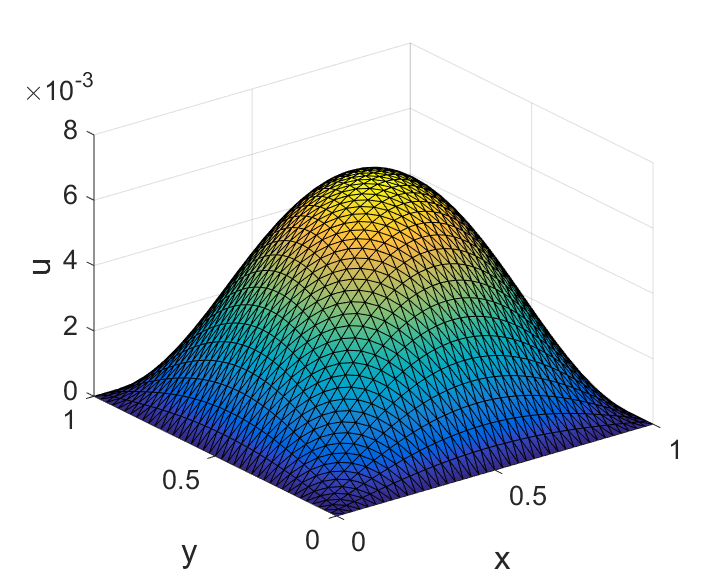}}
\subfigure[{\scriptsize fractional, $t=1.1$}]{
  \includegraphics[trim = .1cm .1cm .1cm .1cm, clip=true,width=0.3\textwidth]{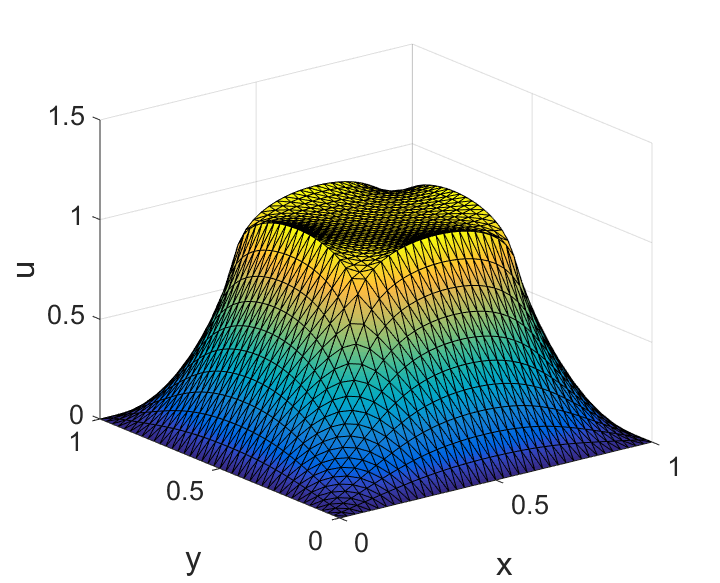}}
\subfigure[{\scriptsize local, $t=1.1$}]{
  \includegraphics[trim = .1cm .1cm .1cm .1cm, clip=true,width=0.3\textwidth]{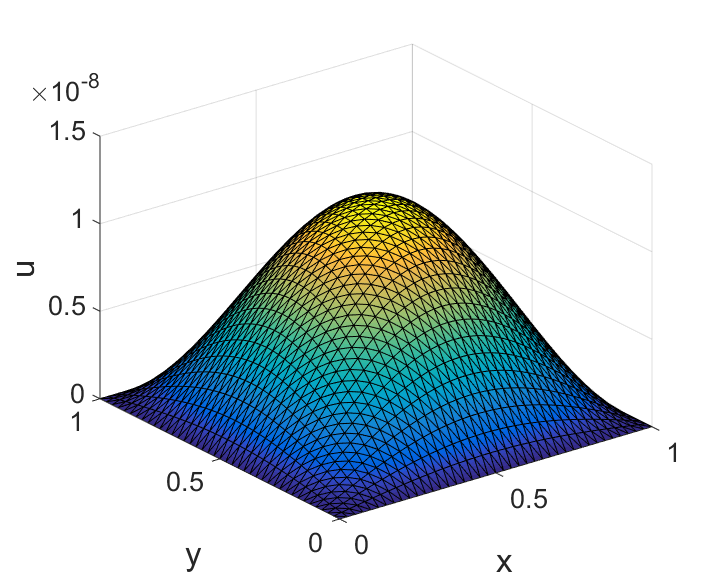}}
\vspace{-0.1cm}
 \caption{Numerical solutions of local, fractional and nonlocal models
 with $\Omega=(0,1)^2$, $\delta=1$ and $\alpha=0.5$ at different time.}\label{fig:case2.2}
\end{figure}

Next, we show a two-dimensional example with $\Om=(0,1)^2$ and an initial data given by
the Dirac-delta measure concentrated on the boundary
of a smaller square $[\frac14,\frac34]\times[\frac14,\frac34]$. Comparisons with the
corresponding solutions of the fractional and local diffusion models are all provided.

\section{Additional discussions on nonlocal-in-time models}\label{sec:3}
As we can see from the examples presented in this paper, the nonlocal-in-time
model (\ref{eqn:fde}) is a simple modification of the traditional dynamic systems where the
instant rate of change ($\partial_t u$) is replaced by a nonlocal rate of change
$(\Gd u)$.  The local limit, i.e., the standard derivative, is simply the extreme case
where the kernel function $\rho_\delta(s)$ degenerates into a singular point measure at $s=0$.
On the other hand, we have also mentioned that by picking suitable fractional
type kernels with $\delta\rightarrow \infty$, the nonlocal-in-time model can recover
fractional differential equations as well.
We thus advocate the equation (\ref{eqn:fde}) as a more general model, whose mathematical theory (such as
existence, uniqueness and regularity of solutions), along with
numerical analysis of some discrete approximations,
can be found in \cite{DuYangZhou:2016}.

There have been a large number of studies on time-fractional dynamics, which are used to
model the continuous time random walk of particles in heterogeneous media \cite{Metzler:2000,Berkowitz:2002,metzler2004restaurant}.
The fractional derivative arising in fractional sub-diffusion model results from a fractional power-law type waiting time probability.
However, such specific choices of memory kernels are rather restrictive. One may argue that perhaps
there remains a lack of compelling evidence that the nature is  confined by such limited forms of kernels.
For example,  let us consider the situation where
particles may explore some environment that is a kind of labyrinthine (bearing short-time sub-diffusion) but homogenizes
on large length scales (resulting in a long-time normal diffusion).
In recent experimental studies of such systems, there have been a number reports on different diffusion regimes
from short-time sub-diffusion to long-time normal diffusion \cite{jeon2012anomalous, he2016dynamic}.
This motivates us to look for a simple and effective model to capture the transient dynamics.
Naturally, there may be different mathematical models for such crossover phenomena. For example,
besides the recently proposed diffusing diffusivity model (cf. \cite{chubynsky2014diffusing,grebenkov2019unifying,jain2017diffusing}),
 a variable-order fractional diffusion equation
$$ \partial_t^{\alpha(t)} u -\Delta u = 0$$
has been used
in \cite{Sun:2017} with a time-dependent order $\alpha(t)$.
However, analyzing such model (theoretically or numerically) remains a
difficult task from a mathematical perspective, besides the challenge for further physical validation.

In contrast, the nonlocal-in-time model (\ref{eqn:fde}) with a constant finite history dependence leads to
a more direct and intuitive way to interpret and to capture the crossover behavior.
Moreover, although different models
may produce similar behavior on MSD, other statistics and spatial/temporal patterns of the solutions
may differ \cite{metzler2019brownian,sposini2018random}. Thus, it is important to conduct more mathematical
investigations and numerical simulations, like the studies presented here, in order
to provide a deeper understanding of the underlying process.

From a modeling perspective, the nonlocal-in-time model (\ref{eqn:fde}) can also be explained by a model of random walk with trapping. For a
truncated fractional kernel function of the type (\ref{fractype}), 
 the probability intensity of the trapping event
is related to  $c_\alpha s^{-\alpha-1}\chi_{(0,\delta)}=:\rho_\delta(s)/s$.
As an illustration, let us focus on a simple 1-D case.
First, in discrete time steps with a uniform step size $\tau$, a random walker is assumed to stay
in its current position, or move to
one of its nearest neighbor sites with length $h$ and in a random left or right direction.
Let $u(x,t)$ be the probability that the walker appears at position $x$ and time $t$,
and $\eta(x,t)$ be the probability that the  walker has just arrived at $x$ and $t$.
By assuming that the maximal waiting time is $\delta$, we then get
\begin{equation}\label{eqn:trap}
u(x,t) = \sum_{k=1}^{\delta/\tau} \omega_k \eta(x,t-k\tau)+ \eta(x,t ),
\end{equation}
where $\omega_k$ is the probability that the  walker waits for at least $k\tau$
after arriving at a position.
We let $p_k$ be the probability that the  walker stops exactly $k\tau$
after its arrival, which for $k=1,2,...,\delta/\tau$, is assumed to be of the form:
\begin{eqnarray*}
    p_{k-1} = \Big(\frac{1}{k\tau}\int_{(k-1)\tau}^{k\tau} s^{-\alpha}  \,ds \Big)/
   \Big( \sum_{k=1}^{\delta/\tau} \frac{1}{k\tau} \int_{(k-1)\tau}^{k\tau} s^{-\alpha} \,ds \Big).
\end{eqnarray*}
Then, it is obvious that
$$    \omega_k = 1 - \sum_{i=0}^{k-1} p_k.$$
The discrete  model of random  walk with trapping  gives that for $h$ small,
\begin{eqnarray*}
\eta(x,t) &= \sum_{k=1}^{\delta/\tau} \frac12(\eta(x-h,t-k\tau)+\eta(x+h,t-k\tau)) p_{k-1}\\
&  \approx \sum_{k=1}^{\delta/\tau} (\eta (x, t-k\tau) + \frac {h^2} 2 \eta_{xx}(x, t-k\tau)) p_{k-1}.
\end{eqnarray*}

Let $ \Gd^\tau$ be defined by
\begin{eqnarray*}
    \Gd^\tau u(t)&= u(t) \Big(\sum_{k=1}^{\delta/\tau} \frac{1}{k\tau} \int_{(k-1)\tau}^{k\tau} \rho_\delta(s) \,ds \Big)
    -\sum_{k=1}^{\delta/\tau} \frac{u(t-k\tau)}{k\tau} \int_{(k-1)\tau}^{k\tau} \rho_\delta(s) \,ds,
\end{eqnarray*}
with  $\rho_\delta(s)= (1-\alpha) \delta^{\alpha-1} s^{-\alpha} $.
Then a simple calculation yields  that
$   \Gd^\tau u(t)  \approx \Gd u(t) $ for small $\tau$,
and for 
\begin{eqnarray*}
c_{\delta,\alpha} = (1-\alpha) \delta^{\alpha-1} \sum_{k=1}^\infty \frac{1}{k} \int_{k-1}^k s^{-\alpha}\,ds
\end{eqnarray*}
we can get
\begin{equation*}
   \sum_{k=1}^{\delta/\tau} \frac{1}{k\tau} \int_{(k-1)\tau}^{k\tau} \rho_\delta(s) \,ds \approx c_{\delta,\alpha} \tau^{-\alpha}.
\end{equation*}
Therefore, we observe that for
%
for small $\tau^\alpha\approx h^2$,
\begin{equation}\label{eqn:der2}
\Gd^\tau \eta (x,t) - \frac{c_{\delta,\alpha}}{2} \frac{h^2}{\tau^\alpha} \eta_{xx} (x,t)=0.
\end{equation}
The nonlocal-in-time model (\ref{eqn:fde}) also follows from (\ref{eqn:der2})
by applying $\Gd^\tau$ on the equation (\ref{eqn:trap}), using
(\ref{eqn:der2}) and letting $\tau^\alpha\approx h^2\rightarrow0$.

\begin{figure}[hbt!]
\centering
{\includegraphics[trim = .1cm .1cm .1cm .1cm, clip=true,width=0.45\textwidth]{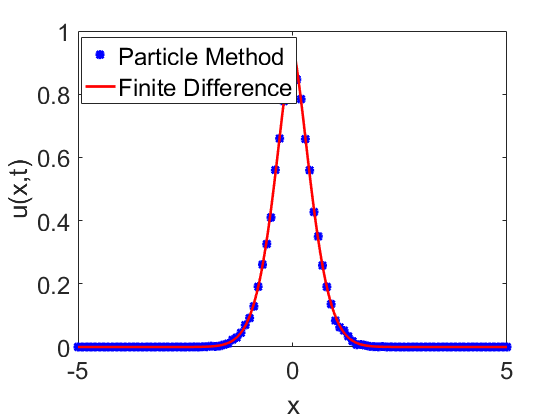}}
{\includegraphics[trim = .1cm .1cm .1cm .1cm, clip=true,width=0.45\textwidth]{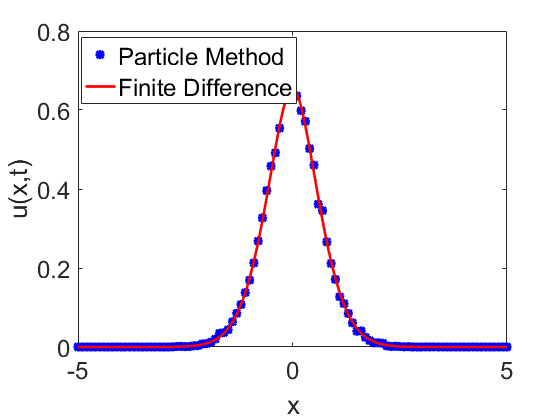}}
{\includegraphics[trim = .1cm .1cm .1cm .1cm, clip=true,width=0.45\textwidth]{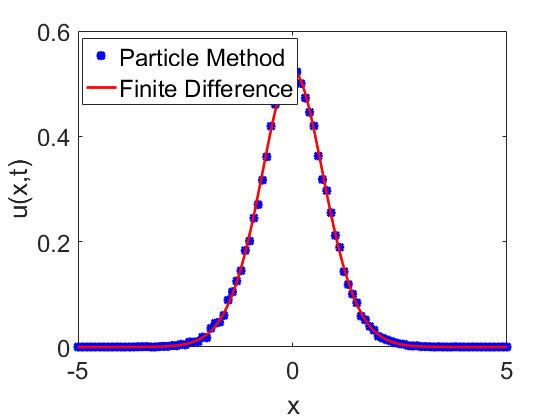}}
{\includegraphics[trim = .1cm .1cm .1cm .1cm, clip=true,width=0.45\textwidth]{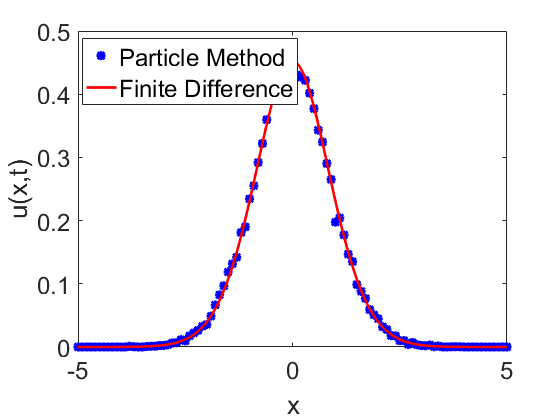}}
\vspace{-0.1cm}
  \caption{Numerical solutions of the nonlocal-in-time diffusion equation
  with $\delta=0.2$, $\alpha=0.75$ and Dirac-Delta initial condition, at $t=0.1$, $0.2$, $0.3$ and $0.4$, respectively.}\label{fig:mc}
\end{figure}

To substantiate this stochastic explanation, we compared numerical solution of (\ref{eqn:der2})
computed by finite different method
with the Monte-Carlo solution using the particle method, in Fig. \ref{fig:mc}.
We consider the infinite domain in one space dimension $\Omega=\mathbb{R}$, and let the initial (historical) data be the
Dirac-delta measure concentrated at $x=0$
i.e., we only track the movement of a single particle, starting from the position $x=0$.
Numerical results indicate that the Monte-Carlo solution fits the finite difference solution very well,
and this supports our theoretical results.

\section{Conclusion}\label{sec:4}
In conclusion, through the computation of solutions and their statistics such as the mean square
displacement, and through the comparisons of them with
the local and fractional counterparts,  this paper shows that the nonlocal-in-time diffusion
with a finite memory can be a very  effective model that
provides an intermediate case between normal and fractional diffusions
and can serve to describe effectively various processes involving crossovers of different
diffusion regimes. The model is showing promises in capturing experimental
observations of diffusion in heterogeneous media without introducing complications and
heterogeneities in the model themselves. The essence lies in the finite memory effect
and how it compares with the overall dynamic history. One may naturally ask
how the memory kernel and the horizon should be chosen, which becomes an interesting
inverse problem to be further investigated. Other interesting issues to be studied
include the connection to other models of anomalous diffusion such as those discussed in
the above and in
section \ref{sec:intro} and
the replacement of normal diffusion operator (the Laplacian $\Delta$) by nonlocal
diffusion operators in the spatial directions. The latter is often associated with
super-diffusion \cite{Du:rev,Defterli:2015},  so that a combined nonlocal in time and space diffusion model
might effectively describe sub-, normal- and super-diffusion regimes.

\bibliographystyle{abbrv}

\end{document}